\documentclass[a4paper,12pt,oneside]{article}
\usepackage[english]{babel}
\usepackage[latin1]{inputenc}
\usepackage{color,a4,pstricks}
\usepackage{graphicx}
\usepackage{dcolumn}
\usepackage{bm}
\usepackage{hyperref} 
\usepackage{mathrsfs} 
\usepackage{amsmath}
\usepackage{amssymb}
\usepackage{amsmath}

\begin{document}
\baselineskip=18pt

\newtheorem{theorem}{Theorem}
\newtheorem{lemma}[theorem]{Lemma}


\def\P{{\bf P}}
\def\E{{\bf E}}
\def\R{{\bf R}}
\def\ee{\mathrm{e}}
\def\d{\, \mathrm{d}}

\begin{center}
{\Large Large deviations for the rightmost position}

\medskip

{\Large in a branching Brownian motion}
\end{center}

\bigskip
\bigskip

\centerline{Bernard Derrida\footnote{\scriptsize Coll\`ege de France, 11 place Marcelin Berthelot, F-75231 Paris Cedex 05, France, and Laboratoire de Physique Statistique, \'Ecole Normale Sup\'erieure, Universit\'e Pierre et Marie Curie, Universit\'e Denis Diderot, CNRS, 24 rue Lhomond, F-75231 Paris Cedex 05, France, {\tt derrida@lps.ens.fr}} and Zhan Shi\footnote{\scriptsize LPMA, Universit\'e Pierre et Marie Curie, 4 place Jussieu, F-75252 Paris Cedex 05, France, {\tt zhan.shi@upmc.fr}}}

\bigskip

\centerline{\it Dedicated to Professor Valentin Konakov}

\centerline{\it on the occasion of his 70th birthday}

\bigskip
\bigskip

{\leftskip=1truecm \rightskip=1truecm \baselineskip=15pt \small

\noindent{\slshape\bfseries Summary.} We study the lower deviation probability of the position of the rightmost particle in a branching Brownian motion and obtain its large deviation function.

\bigskip


\bigskip


} 

\bigskip
\bigskip

\section{Introduction}
\label{s:introduction}

The question of the distribution of the position $X_{\max}(t)$ of the rightmost particle in a branching Brownian motion (BBM) has a long history  in probability theory \cite{McKean,bramson78,bramson83,chauvin-rouault,rouault,HS,zeitouni,berestycki,Shi,Bovier} and in physics \cite{DSpohn,MK,MM,RMS}. 

By branching Brownian motion, we mean that the system starts with a single particle at the origin which performs a Brownian motion with variance $\sigma^2$ at time $1$, and branches at rate $1$ into two independent Brownian motions which themselves branch at rate $1$ independently, and so on. For such a BBM, one knows since the work of McKean \cite{McKean}
that 
$$
u(x, \, t)
:=
\P (X_{\max}(t) \le x),
$$

\noindent satisfies the F-KPP (Fisher--Kolmogorov-Petrovskii-Piskounov) equation
\begin{equation}
{\partial u \over \partial t} = {\sigma^2 \over 2} {\partial^2 u \over \partial x^2} +  u^2 -u
\label{F-KPP}
\end{equation}
with the initial condition $u(x,0)= {\bf 1}_{\{ x\ge 0\}}$.
It is also known since the works of Bramson \cite{bramson78,bramson83} that in the long time limit
\begin{equation}
u(x+m(t)\sigma, \, t) \to F(x)
\label{Bramson}
\end{equation} where $F(z) $
is a traveling wave solution of 
$$ {\sigma^2 \over 2} F''  + \sqrt{2\sigma^2\,} F' + F^2 - F=0$$
and 
\begin{equation}
m(t) := \sqrt{2}\, t - {3\over 2 \sqrt{2} } \ln t\, .
\label{m}
\end{equation}

\noindent This implies in particular that
$$
\lim_{t\to \infty} \frac{X_{\max}(t)}{t}
=
\sqrt{2\sigma^2} \, ,
\qquad\hbox{\rm in probability.}
$$

\noindent [The convergence also holds almost surely.]

In 1988 Chauvin and Rouault \cite{chauvin-rouault,rouault} proved  a large deviation result for $X_{\max}(t) / t > \sqrt{2\sigma^2}$, namely, that for $v > \sqrt{2\sigma^2}$
\begin{equation}
    \ln \left[ \P\left({X_{\max}(t)  \over t} > v \right)\right] \sim  t \left( 1 -{v^2 \over 2 \sigma^2} \right) . 
    \label{Chauvin-Rouault}
\end{equation}

\noindent In (\ref{Chauvin-Rouault}) and everywhere below, the symbol $\sim$ means that
\begin{equation}
    \label{symbol}
    \lim_{t \to \infty} \frac{\ln \P( X_{\max}(t) >vt)}{t ( 1 -{v^2 \over 2 \sigma^2} )} 
    = 
    1 \ .  
\end{equation}

\noindent Here  we are interested in the {\it lower deviation} probability $\P( X_{\max}(t) \le vt)$ for each $v\in (-\infty, \, \sqrt{2\sigma^2}\, )$. It turns out that $v/\sqrt{2\sigma^2}$ is an important parameter, so we fix $\alpha\in (-\infty, \, 1)$, and study
$$
\P( X_{\max}(t) \le \alpha\sqrt{2\sigma^2}\, t) ,
$$

\noindent when $t\to \infty$. 

Throughout the paper, we write
\begin{equation}
    \rho := \sqrt{2}-1 \, .
    \label{rho}
\end{equation}

\medskip

\begin{theorem}
\label{t:BBM}

Let $X_{\max}(t)$ denote the rightmost position of the BBM at time $t$. Then for all $\alpha\in (-\infty, \, 1)$,
\begin{equation}
    \ln \P( X_{\max}(t) \le \alpha\sqrt{2\sigma^2}\, t)
    \sim
    - t \, \psi (\alpha)
    \label{BBM}
\end{equation}

\noindent where
\begin{equation}
    \psi (\alpha)
    =
    \begin{cases}
      2\rho(1- \alpha)\, , 
       &\hbox{ if } \alpha\in [-\rho, \, 1) \, ,
       \\
      1+ \alpha^2 \, ,
       &\hbox{ if } \alpha\in (-\infty, \, -\rho] \, .
    \end{cases}
    \label{Res-NBBM}
\end{equation}

\end{theorem}

\medskip

\begin{figure}[h]
\centerline{\includegraphics[scale=.3]{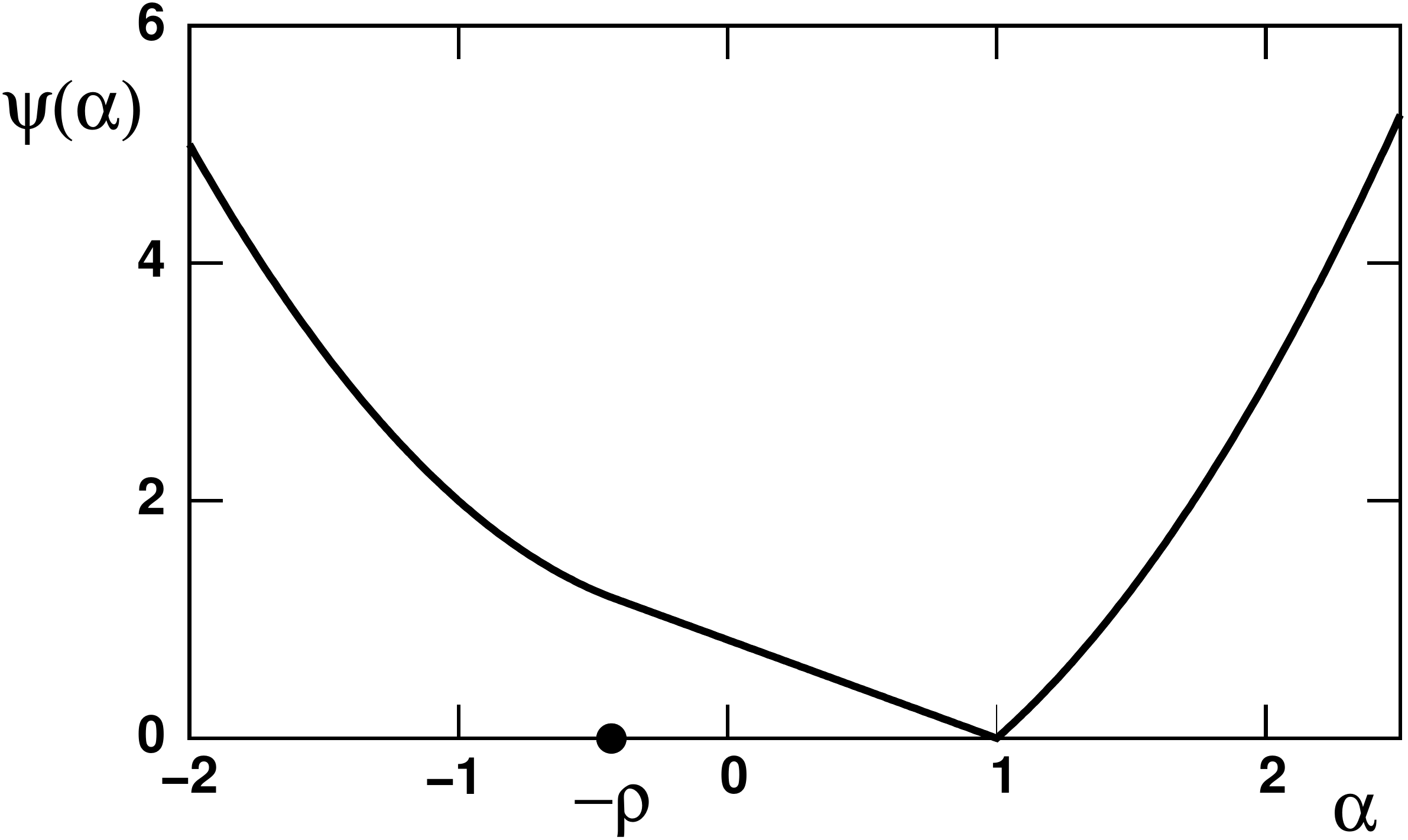}}
\caption{\small The large deviation function of the position of the rightmost particle of a branching Brownian motion.  
The expression of $\psi(\alpha)$ is non-analytic at $\alpha=-\rho=1-\sqrt{2}$ and at $\alpha=1$.
 }
\label{ld-figure}
\end{figure}

Together with Theorem \ref{t:BBM} and the upper large deviation probability in \eqref{Chauvin-Rouault}, a routine argument (proof of Theorem III.3.4 in den Hollander~\cite{denhollander}, proof of Theorem 2.2.3 in Dembo and Zeitouni~\cite{dembo-zeitouni}) yields the following formalism of large deviation principle: the family of the distributions of $\frac{X_{\max}(t)}{\sqrt{2\sigma^2}\, t}$, for $t\ge 1$, satisfies the large deviation principle on $\R$, with speed $t$ and with the rate function $\psi(\alpha)$ (shown in   Figure \ref{ld-figure})
\begin{equation}
    \psi (\alpha)
    =
    \begin{cases}
      1+\alpha^2 \, , 
       &\hbox{ \rm if } \alpha \le -\rho \, ,
       \\
      2\rho(1- \alpha) \, ,
       &\hbox{ \rm if } -\rho \le \alpha \le 1 \, ,
       \\
      \alpha^2 -1 \, ,
       &\hbox{ \rm if } \alpha \ge 1 \, ,
    \end{cases}
    \label{psi}
\end{equation}

\noindent i.e., for any closed set $F \subset \R$ and open set $G \subset \R$,
\begin{eqnarray*}
    \limsup_{t\to \infty} \frac1t \ln \P \Big( \frac{X_{\max}(t)}{\sqrt{2\sigma^2}\, t} \in F \Big) 
 &\le& - \inf_{\alpha \in F} \psi (\alpha),
    \\
    \liminf_{t\to \infty} \frac1t \ln \P \Big( \frac{X_{\max}(t)}{\sqrt{2\sigma^2}\, t} \in G \Big) 
 &\ge& - \inf_{\alpha \in G} \psi (\alpha) .
\end{eqnarray*}

Let us also mention that Proposition 2.5 of Chen~\cite{chen} (recalled as Lemma \ref{l:chen} in Section \ref{s:ub} below) implies that for all $\alpha<1$,
$$
\psi (\alpha)
\ge
\frac{1-\alpha}{6} \, ,
$$
\noindent which is in agreement with \eqref{Res-NBBM}.

The reason for the non-analyticity of  $\psi(\alpha)$ in  (\ref{psi}) at $\alpha= - \rho$ is that, as we will see it in sections 2 and 3, for $\alpha < - \rho$ the events which  dominate are those where the  initial particle does not branch or branches at a very late time (at a time $\tau$ very close to $t$) while in the range $-\rho < \alpha < 1$ the first branching event occurs at a time $\tau \sim (1-\alpha) t /\sqrt{2}$.

The rest of the paper is as follows. Sections \ref{s:lb} and \ref{s:ub} are devoted to the proof of the lower bound and the upper bound, respectively, for the probability in Theorem \ref{t:BBM}. In Section \ref{s:conclusion_and_remarks}, we present some further remarks.

\section{Lower bound}
\label{s:lb}

Fix $v\in (-\infty, \, \sqrt{2\sigma^2}\,)$. We prove the lower bound in the deviation probability, by considering a special event described as follows: The initial particle does not produce any offspring during time interval $[0, \, \tau]$ and is positioned at $y\in (-\infty, \, vt-\sqrt{2\sigma^2}\, (t-\tau)-1]$ at time $\tau$; then, at time $t$, the maximal position lies in $(-\infty, \, vt)$. As such, we get
\begin{eqnarray}
 &&\P(X_{\max}(t) \le vt)
    \nonumber
    \\
 &\ge& \ee^{-\tau} \, \int_{-\infty}^{vt-\sqrt{2\sigma^2}\, (t-\tau)-1} \frac{\d y}{\sqrt{2\pi \sigma^2 \tau}} \, \ee^{- \frac{y^2}{2\sigma^2 \tau}} \, \P (X_{\max}(t-\tau) < vt - y) \, .
    \label{pf:lb_eq1}
\end{eqnarray}

\noindent Note that for $y\in (-\infty, \, vt-\sqrt{2\sigma^2}\, (t-\tau)-1]$, we have $vt-y \ge \sqrt{2\sigma^2}\, (t-\tau)+1$, so 
$$
\P (X_{\max}(t-\tau) \le vt - y)
\ge
\P (X_{\max}(t-\tau) \le \sqrt{2\sigma^2}\, (t-\tau) +1) \, .
$$

\noindent Let $m(t) := \sqrt{2}\, t - {3\over 2 \sqrt{2} } \ln t$ be as in \eqref{m}. By \eqref{Bramson}, for any $z\in \R$, $\P (X_{\max}(s) \le m(s)\sigma + z)$ converges, as $s\to \infty$, to a positive limit (which depends on $z$). This yields the existence of a constant $c>0$ such that 
$$
\P (X_{\max}(t-\tau) \le \sqrt{2\sigma^2}\, (t-\tau) +1) 
\ge 
c,
$$

\noindent for all $\tau \in [0, \, t]$. [The presence of $+1$ in $\sqrt{2\sigma^2}\, (t-\tau) +1$ is only to ensure the positivity of the probability when $\tau$ equals $t$ or is  very close to $t$.] Going back to \eqref{pf:lb_eq1}, we get that for all $\tau \in (0, \, t]$,
$$
\P(X_{\max}(t) \le vt)
\ge
c\, \ee^{-\tau} \, \int_{-\infty}^{vt-\sqrt{2\sigma^2}\, (t-\tau)-1} \frac{1}{\sqrt{2\pi \sigma^2 \tau}} \, \ee^{- \frac{y^2}{2\sigma^2 \tau}} \d y \, .
$$

\noindent Hence
\begin{equation}
    \P(X_{\max}(t) \le vt)
    \ge
    c\sup_{\tau\in (0, \, t]} \Big\{ \ee^{-\tau} \, \int_{-\infty}^{vt-\sqrt{2\sigma^2}\, (t-\tau)-1} \frac{1}{\sqrt{2\pi \sigma^2 \tau}} \, \ee^{- \frac{y^2}{2\sigma^2 \tau}} \d y \Big\} \, .
    \label{pf:lb_eq2}
\end{equation}

\noindent We now use the following result.

\begin{lemma}
\label{l:calcul_integrale}

 For $v< \sqrt{2\sigma^2}$ and $t\to \infty$,
 $$
 \ln \Big( \sup_{\tau\in (0, \, t]} \Big\{ \ee^{-\tau} \, \int_{-\infty}^{vt-\sqrt{2\sigma^2}\, (t-\tau)-1} \frac{1}{\sqrt{2\pi \sigma^2 \tau}} \, \ee^{- \frac{y^2}{2\sigma^2 \tau}} \d y \Big\} \Big)
 \sim 
 - \varphi(v) t ,
 $$
 where
 \begin{equation}
     \varphi(v)
     :=
     \begin{cases}
       2\rho(1-\alpha) , &\hbox{ \rm if } \alpha\ge -\rho,
       \\
       1+\alpha^2, &\hbox{ \rm if } \alpha\le -\rho \, ,
     \end{cases}
     \label{phi}
 \end{equation}
 with $\alpha:= \frac{v}{\sqrt{2\sigma^2}} <1$ and $\rho := \sqrt{2}-1$ as before.

\end{lemma}

The proof of Lemma \ref{l:calcul_integrale} is quite elementary (as $\ln (\int_{-\infty}^z \ee^{-y^2} \d y) \sim - z^2$ for $z\to -\infty$, and $\int_{-\infty}^z \ee^{-y^2} \d y$ is greater than a positive constant if $z\ge 0$). We only indicate the optimal value of $\tau$:
\begin{equation}
    \tau
    =
    \begin{cases}
     \frac{1-\alpha}{\sqrt{2}}\, t +o(t), 
      &\hbox{ \rm if } \alpha\ge -\rho ,
       \\
     t +o(t),
      &\hbox{ \rm if } \alpha\le -\rho\, .
    \end{cases}
    \label{tau:valeur_optimale}
\end{equation}

By \eqref{pf:lb_eq2} and Lemma \ref{l:calcul_integrale}, we obtain:
$$
\liminf_{t\to \infty} \frac1t \, \ln \P(X_{\max}(t) \le vt)
\ge
- \varphi(v) \, ,
$$

\noindent with $\varphi(v)$ as in \eqref{phi}. This yields the desired lower bound for the probability in the theorem, as $\varphi(v)$ coincides with $\psi(\alpha)$ defined in \eqref{psi}.

\section{Upper bound}
\label{s:ub}

We now look for the upper bound in the deviation probability. 
Fix $x=vt$ with $v<\sqrt{2\sigma^2}\,$. Let
$$
u(x, \, t)
:=
\P (X_{\max}(t) \le x),
$$

\noindent as before. Considering the event that the first branching time is $\tau$, we have
\begin{eqnarray*}
    u(x, \, t)
 &=& \int_{-\infty}^x \frac{\d y}{\sqrt{2\pi \sigma^2 t}} \, \ee^{-t-\frac{y^2}{2\sigma^2 t}}
    \\
 && \quad 
    +
    \int_0^t \d \tau \int_{-\infty}^\infty \frac{\d y}{\sqrt{2\pi \sigma^2 \tau}} \, \ee^{-\tau-\frac{y^2}{2\sigma^2 \tau}}\, u^2 (x-y, \, t-\tau)\, ,
\end{eqnarray*}

\noindent the first term on the right-hand side originating from the event that the first branching time is greater than $t$. [It is easy to check that this expression satisfies the F-KPP equation \eqref{F-KPP}.] We also have a lower bound for $u(x, \, t)$ by considering only the event that there is no branching up to time $\tau$: For any $\tau\in [0, \, t]$, 
$$
    u(x, \, t)
    \ge
    \int_{-\infty}^\infty \frac{\d y}{\sqrt{2\pi \sigma^2 \tau}} \, \ee^{-\tau-\frac{y^2}{2\sigma^2 \tau}}\, u (x-y, \, t-\tau)\, .
$$

\noindent Writing $(B_s, \, s\ge 0)$ for a standard Brownian motion (with variance of $B_1$ being $1$), the last two displayed formulas can be expressed as follows:
\begin{eqnarray}
    u(x, \, t)
 &=& \ee^{-t} \, \P( \sigma B_t \le x)
    +
    \int_0^t \ee^{-\tau} \, \E[ u^2(x-\sigma B_\tau, \, t-\tau)] \d \tau \, ,
    \label{pf:ub_eq1}
    \\
    u(x, \, t)
 &\ge& \ee^{-\tau} \, \E[ u(x-\sigma B_\tau, \, t-\tau)],
    \qquad
    \forall \tau\in [0, \, t]\, .
    \label{pf:ub_eq2}
\end{eqnarray}

Consider, for $\tau \in [0, \, t]$,
$$
\Phi(\tau)
:=
\ee^{-\tau} \, \E[ u^2(x-\sigma B_\tau, \, t-\tau)] \, .
$$

\noindent Since $\Phi$ is a continuous function on $[0, \, t]$, there exists $\tau_0 = \tau_0(t, \, x)$ such that
$$
\Phi(\tau_0)
=
\sup_{\tau\in [0, \, t]} \Phi(\tau) \, .
$$

\noindent On the other hand, since $u(\, \cdot \, , 0) = {\bf 1}_{[0, \, \infty)}(\, \cdot \,)$, we have $\ee^{-t} \, \P( \sigma B_t \le x) = \Phi(t)$. So \eqref{pf:ub_eq1} becomes $u(x, \, t) = \Phi(t) + \int_0^t \Phi(\tau) \d \tau$, which is bounded by $(t+1)\sup_{\tau\in [0, \, t]} \Phi(\tau)$. Taking $\tau = \tau_0$ in \eqref{pf:ub_eq2}, it follows from \eqref{pf:ub_eq2} and \eqref{pf:ub_eq1} that
$$
    \ee^{-\tau_0} \, \E[ u(x-\sigma B_{\tau_0}, \, t-\tau_0)]
    \le 
    u(x, \, t)
    \le
    (t+1) \Phi(\tau_0)\, ,
$$

\noindent which can be represented as
\begin{equation}
    \ee^{-\tau_0} \, \E(Y)
    \le
    u(x, \, t)
    \le
    (t+1) \ee^{-\tau_0} \, \E(Y^2) \, ,
    \label{pf:ub_eq3}
\end{equation}

\noindent where
$$
Y
=
Y(x, \, t, \, \sigma)
:=
u(x-\sigma B_{\tau_0}, \, t-\tau_0) \, .
$$

Let us have a closer look at $\E(Y)$. We write
$$
\ee^{-\tau_0} \, \E(Y)
=
A_1
+
A_2 \, ,
$$

\noindent with
\begin{eqnarray*}
    A_1 = A_1(x, \, t, \, \sigma)
 &:=& \ee^{-\tau_0} \, \E[ Y\, {\bf 1}_{\{ Y < \frac{1}{2(t+1)} \}} ] \, ,
    \\
    A_2 = A_2(x, \, t, \, \sigma)
 &:=& \ee^{-\tau_0} \, \E[ Y\, {\bf 1}_{\{ Y \ge \frac{1}{2(t+1)} \}} ] \, .
\end{eqnarray*}

\noindent Then
\begin{eqnarray*}
 &&(t+1) \ee^{-\tau_0}\, \E(Y^2)
    \\
 &=& (t+1) \ee^{-\tau_0}\, \E[ Y^2\, {\bf 1}_{\{ Y < \frac{1}{2(t+1)} \}} ]
    +
    (t+1) \ee^{-\tau_0}\, \E[ Y^2\, {\bf 1}_{\{ Y \ge \frac{1}{2(t+1)} \}} ]
    \\
 &\le& \frac12 \, \ee^{-\tau_0}\, \E[ Y\, {\bf 1}_{\{ Y < \frac{1}{2(t+1)} \}} ]
    +
    (t+1) \ee^{-\tau_0}\, \E[ Y \, {\bf 1}_{\{ Y \ge \frac{1}{2(t+1)} \}} ]\, ,
\end{eqnarray*}

\noindent where, on the right-hand side, we have used the trivial inequality $Y^2 \le Y$ when dealing with the event $\{ Y \ge \frac{1}{2(t+1)} \}$. In other words,
$$
(t+1) \ee^{-\tau_0}\, \E(Y^2)
\le
\frac12 \, A_1
+
(t+1) A_2 \, .
$$

\noindent So by \eqref{pf:ub_eq3}, we obtain
$$
A_1+A_2
\le 
u(x, \, t)
\le
(t+1) \ee^{-\tau_0}\, \E(Y^2)
\le
\frac12 \, A_1
+
(t+1) A_2 \, .
$$

\noindent In particular, this implies $A_1 \le 2t A_2$. As a consequence,
\begin{equation}
    A_2
    \le 
    u(x, \, t)
    \le
    (2t+1) A_2 \, .
    \label{pf:ub_eq4}
\end{equation}

\noindent This yields that $A_2$ has the same asymptotic behaviour as $u(x, \, t)$, as far as large deviation functions are concerned. 

We now look for an upper bound for $A_2$, which, multiplied by $2t+1$, will be served as an upper bound for $u(x, \, t)$. Let us recall the following estimate:

\begin{lemma}
\label{l:chen}

 {\bf (Chen~\cite{chen}, Proposition 2.5)}
 Let $m(t) := \sqrt{2}\, t - \frac{3}{2 \sqrt{2}} \, \ln t$ as in \eqref{m}. There exist two constants $c_1>0$ and $c_2>0$ independent of $\sigma$, such that 
 $$
 \P(X_{\max}(r) \le \sigma m(r) - \sigma z \hbox{ \it for some } r\le \ee^z) 
 \le 
 c_1 \, \ee^{-c_2 z} ,
 $$
 for all sufficiently large $z$. Moreover, one can take $c_2 = \frac{1}{6\sqrt{2}}$.

\end{lemma}

We apply the lemma to $z:= t^{1/3}$, to see that when $t$ is sufficiently large (say $t\ge t_0$), for any $\tau\in [0, \, t]$,
$$
y < \sqrt{2\sigma^2} \, \tau - t^{1/2} 
\; \Rightarrow \;
u(y, \, \tau) < \frac{1}{2(t+1)} \, .
$$

\noindent 
As such, for $t\ge t_0$, we have 
$$
A_2
=
\ee^{-\tau_0} \, \E[ Y\, {\bf 1}_{\{ Y \ge \frac{1}{2(t+1)} \}} ] 
\le
\ee^{-\tau_0} \, \E[ Y\, {\bf 1}_{\{ x-\sigma B_{\tau_0} \ge \sqrt{2\sigma^2} \, (t-\tau_0) - t^{1/2} \}} ] \, .
$$

\noindent Since $Y\le 1$, this yields, for $t\ge t_0$, 
\begin{eqnarray*}
    A_2
 &\le& \ee^{-\tau_0} \, \P (x-\sigma B_{\tau_0} \ge \sqrt{2\sigma^2} \, (t-\tau_0) - t^{1/2})
    \\
 &\le& \sup_{\tau \in [0, \, t]} \Big\{ \ee^{-\tau} \, \P (x-\sigma B_\tau \ge \sqrt{2\sigma^2} \, (t-\tau) - t^{1/2}) \Big\}
    \\
 &=&
    \sup_{\tau \in (0, \, t]} \Big\{ \int_{-\infty}^{x-\sqrt{2\sigma^2} \, (t-\tau) + t^{1/2}} \frac{1}{\sqrt{2\pi \sigma^2 \tau}} \, \ee^{-\tau -\frac{y^2}{2\sigma^2 \tau}} \d y\Big\} \, ,
\end{eqnarray*}

\noindent By \eqref{pf:ub_eq4}, we have therefore, for all sufficiently large $t$,
$$
u(x, \, t)
\le
(2t+1) \sup_{\tau \in (0, \, t]} \Big\{ \int_{-\infty}^{x-\sqrt{2\sigma^2} \, (t-\tau) + t^{1/2}} \frac{1}{\sqrt{2\pi \sigma^2 \tau}} \, \ee^{-\tau -\frac{y^2}{2\sigma^2 \tau}} \d y\Big\} \, .
$$

Recall that $x=vt$. The supremum on the right-hand side has already been estimated in Lemma \ref{l:calcul_integrale} in Section \ref{s:lb}: For $v< \sqrt{2\sigma^2}$ and $t\to \infty$,
$$
\ln \Big( \sup_{\tau \in (0, \, t]} \Big\{ \int_{-\infty}^{x-\sqrt{2\sigma^2} \, (t-\tau) + t^{1/2}} \frac{1}{\sqrt{2\pi \sigma^2 \tau}} \, \ee^{-\tau -\frac{y^2}{2\sigma^2 \tau}} \d y\Big\} \Big)
\sim 
-\varphi(v) t ,
$$

\noindent where $\varphi(v)$ is defined in \eqref{phi}. Note that we have $t^{1/2}$ here (in $x-\sqrt{2\sigma^2} \, (t-\tau) + t^{1/2}$) instead of $-1$ in the lemma; this makes in practice no difference because $t^{1/2} \le \varepsilon t$ (for any $\varepsilon>0$ and all sufficiently large $t$) and we can use the continuity of the function $v\mapsto \varphi(v)$. Consequently, for $x= vt$ with $v< \sqrt{2\sigma^2}$,
$$
\limsup_{t\to \infty} \frac{\ln u(x, \, t)}{t}
\le
-\varphi(v) ,
$$

\noindent which yields the upper bound for the probability in the theorem because $\varphi(v)$ coincides with $\psi(\alpha)$ given in \eqref{Res-NBBM}.

\section{Conclusion and remarks}
\label{s:conclusion_and_remarks}

The main result stated in
(\ref{Res-NBBM}) and (\ref{psi})
of the present  work is the expression of the (lower) large deviation function  $\psi(\alpha)$ of the position of  the rightmost particle of a branching brownian motion.
One remarkable feature of this large deviation function is its non-analyticty
at some particular values $\alpha=-\sqrt{2} +1$ and $\alpha=1$ due to a change of scenario of the dominant contribution to the large  deviation function:
for $\alpha  
 < -\sqrt{2}  
 +1$, 
the dominant event is a single Brownian particle which does not branch up to time $t$; for $ -\sqrt{2} +1 < \alpha < 1$, 
it corresponds to a particle which moves  to position $ -(\sqrt{2} -1) (1- \alpha) \sigma t$ without branching  up to a
 time $t(1-\alpha)/\sqrt{2}$, and then behaves like a normal BBM up to time $t$;
for $\alpha > 1$,  the tree branches normally but one branch moves at the speed $\alpha\sqrt{2\sigma^2}$, faster than the normal speed $\sqrt{2\sigma^2}$. 

Using more {\it heuristic} arguments as in \cite{DMS}, it is possible to determine the time dependence of  the prefactor, for example by  showing \cite{DS} that for  $-\rho < \alpha < 1 $, there exists a constant $c\in (0, \, \infty)$ such that
\begin{equation}
\P( X_{\max}(t) \le \alpha\sqrt{2\sigma^2}\, t)  \sim c\, t^{3(\sqrt{2} -1) \over 2} \ee^{-  \psi(\alpha) t} \ .
\label{conjecture}
\end{equation}
The result of the present work can also be easily extended to more general branching Brownian  motions, where one includes the possibility that a particle branches into more than two particles (for example one could consider that a particle branches into $k$ particles with probablity $p_k$). It can also be extended to branching random walks. In all these cases, one finds \cite{DS}  as in 
(\ref{Res-NBBM}) and (\ref{psi}) three different regimes with the same scenarios as described above.

It is however important to notice that  expressions 
(\ref{Res-NBBM}) and (\ref{psi})
  of the large deviation function $\psi(\alpha)$ for $\alpha <1$ depend crucially on the fact that one starts initially with a single particle and that branchings occur at random times according to Poisson processes. If instead one starts at time $t=0$ with several particles in \cite{meerson-sasorov} or if the distribution of the branching times is not exponential (for example in the case of a branching random walk generated by a regular binary tree where at each (integer) time step each particle branches into two particles),  $\P (X_{\max} \le v t)$ might decay faster than an exponential of time.

Recently there has been a renewed interest in the understanding of the extremal  process 
and  in particular of the measure seen at the tip of the branching Brownian motion \cite{LS,BD1,BD2,ABBS,ABK,SK}. We think that it would be interesting to investigate how this extremal process is modified when it is conditioned on the position of the rightmost particle, i.e., how it depends on the parameter $\alpha$.

\end{document}